\documentclass[11pt]{amsart}
\usepackage{cite,enumerate,setspace,comment}
\usepackage{MnSymbol}
\usepackage{paralist}
\usepackage{amsmath}
\usepackage{verbatim}
\usepackage{esint}
\usepackage{enumitem}
\theoremstyle{plain}
\newtheorem{theorem}{Theorem}[section]

\newtheorem{cor}[theorem]{Corollary}

\theoremstyle{definition}

\theoremstyle{remark}

\numberwithin{equation}{section}

\allowdisplaybreaks[1]
\begin{document}

\title [A.e. convergence of Fourier series] {Almost everywhere convergence of Fourier series on compact connected Lie groups}

\author{David Grow and Donnie Myers}
\address {Missouri University of Science and Technology\\ Department of Mathematics and Statistics\\ Rolla, Missouri 65409-0020, USA.}

\let\oldthefootnote\thefootnote
\renewcommand{\thefootnote}{\fnsymbol{footnote}}
\footnotetext{Email addresses: {grow@mst.edu}, {dfm@mst.edu}}
\let\thefootnote\oldthefootnote

\begin{abstract}
We consider the open problem: Does every square-integrable function $f$ on a compact, connected Lie group
$G$ have an almost everywhere convergent Fourier series?  We prove a general theorem from which it follows
that if the integral modulus of continuity of $f$ is $O(t^{\alpha})$ for some $\alpha > 0$ then
the Fourier series of $f$ converges almost everywhere on $G$.  In particular, the Fourier series of any
$\alpha$-H\"{o}lder continuous function on $G$ converges almost everywhere. On the other hand, we show that
to each countable subset $E$ of $G = SU(2)$ and each $0 < \alpha < 1$ there corresponds an $\alpha$-H\"{o}lder
continuous function on $SU(2)$ whose Fourier series diverges on $E$.
\end{abstract}



\maketitle

\section{Introduction}

The Peter-Weyl theorem suggests the study of the formal Fourier series \\${\sum d_{\lambda}(\chi_{\lambda}\ast f)}$ of a
square-integrable function $f$ on a compact, connected Lie group $G$.  Here the sum is over the equivalence classes
of continuous irreducible unitary representations of $G$, $d_{\lambda}$ is the degree of the representation, and
$\chi_{\lambda}$ is its character.  The vast literature of Fourier analysis on $G$ is primarily concerned with mean
convergence or divergence of the Fourier series of $f \in L^p(G)$ (e.g.\cite{GST,GLZ,Mi,S,ST}), uniform or absolute
convergence of the partial sums if $f$ is smooth (e.g.\cite{Cl,GLZ,Ma1,DFM,MG,R1,R2,Sug,T}), almost everywhere convergence
or divergence of the partial sums if $f$ is a central function in $L^p(G)$ (e.g.\cite{CGTV,GT,MN}),
and uniform, mean, or almost everywhere summability of the partial sums if $f$ belongs to various subspaces of $L^1(G)$
(e.g.\cite{CT,CF,CGT,GLZ,Z}). The aim of this work is to advance the study of almost everywhere convergence or divergence of
Fourier partial sums of nonsmooth, possibly noncentral functions in $L^2(G)$.

Relying on Jackson's theorem for compact, connected Lie groups \cite{CK} and a general version of the Rademacher-Menshov
theorem \cite{M}, we show that if $f$ in $L^2(G)$ has an integral modulus of continuity $\Omega(f,\cdot)$ satisfying
\begin{equation*}
\int_0^1\frac{\Omega^2(f,t)}{t}dt < \infty,
\end{equation*}
then the sequences of polyhedral Fourier partial sums $\{S_{N}f(x)\}_{N=1}^{\infty}$ and spherical Fourier partial sums
$\{\widetilde{S}_{N}f(x)\}_{N=1}^{\infty}$ converge to $f(x)$ almost everywhere on $G$.  In particular, if $f$ is an
$\alpha$-H\"{o}lder continuous function on $G$ for some $\alpha > 0$, or more generally if $\Omega(f,t) = O(t^{\alpha})$
for some $\alpha >0 $, then these Fourier partial sums of $f$ converge to $f(x)$ almost everywhere on $G$.

On the other hand, consider the two-dimensional special unitary group $G = SU(2)$.  We show that to each $\alpha$ in $(0,1)$
and each countable subset $E$ of $SU(2)$ there corresponds an $\alpha$-H\"{o}lder continuous function on $SU(2)$ whose
Fourier partial sums diverge at each $x$ in $E$.  Since it is possible to arrange that such a set $E$ is dense in $SU(2)$,
the Fourier partial sums of the corresponding function are divergent at infinitely many points in every nonempty open subset
of $SU(2)$, despite the fact that the Fourier partial sums of such a function converge almost everywhere on $SU(2)$.  It is
worth noting in the case $\alpha = 1$, i.e. when $f$ is Lipschitz continuous on $SU(2)$, the Fourier partial sums of
$f$ converge uniformly to $f$ on $SU(2)$ \cite{DFM,MG}.

It is an open problem whether
\begin{equation}\label{aeconvergence}
  \underset{N\rightarrow\infty}{\lim}S_{N}f(x) = f(x)
\end{equation}
holds almost everywhere for every $f$ in $L^2(G)$. A general theorem of Stanton and Tomas \cite{ST} for compact,
connected, semisimple Lie groups $G$ shows that if $p < 2$ then there correspond an $f$ in $L^p(G)$ and a subset
$E$ of $G$ of full measure such that (\ref{aeconvergence}) fails for all $x$ in $E$.  However, Carleson's
celebrated proof of Lusin's conjecture \cite{Ca} guarantees (\ref{aeconvergence}) holds almost everywhere when $f$ is
square-integrable on the circle group $\mathbb{T}$, and this was extended to $\mathbb{T}^n$ for $n \geq 2$ \cite{F1,Sj,Te}.
Furthermore, it follows from a result of Pollard \cite{P} on Jacobi series that if $f$ is a central function in $L^p(SU(2))$
for some $p > 4/3$, then (\ref{aeconvergence}) holds almost everywhere. Finally, the Peter-Weyl theorem implies that to
every compact, connected Lie group $G$ and each $f$ in $L^2(G)$ there corresponds an increasing sequence $\{N_j\}$ of
positive integers such that
\begin{equation*}
  \underset{j\rightarrow\infty}{\lim}S_{N_j}f(x) = f(x)
\end{equation*}
for almost every $x$ in $G$. This lends some hope for a positive answer to the almost everywhere convergence problem
(\ref{aeconvergence}), as do the results in this paper.

\section{Notation and Preliminaries}

We use primarily the notation of \cite{BD}. Let $G$ be a compact, connected Lie group with Haar measure $\mu$, let
$T$ be a fixed maximal torus in $G$, and let $\mathfrak{g}$ and $\mathfrak{t}$ be the Lie algebras of $G$ and $T$,
respectively.  Choose an inner product $\langle \cdot,\cdot \rangle$ on $\mathfrak{g}$ which is invariant under the
adjoint action of $G$ on $\mathfrak{g}$.  (When $G$ is semisimple we may choose $\langle \cdot,\cdot \rangle$ to be
$-B(\cdot,\cdot)$ where $B$ is the Killing form.)  This provides inner products on $\mathfrak{t}$ and $\mathfrak{t}^{*}$
which are invariant under the Weyl group $W = N(T)/T$.  Let $I = \{H \in \mathfrak{t}: \exp(H) = 1\}$ be the integral
lattice, let $I^{*} = \{\lambda \in \mathfrak{t}^{*}: \lambda(H) \in \mathbb{Z} \text{ for all } H \in I\}$ be the lattice
of integral forms, and let $P \subset I^{*}$ be the set of real roots of $G$ with basis $\{\alpha_1,...,\alpha_l\}$,
$P_{+}$ the set of positive roots, and $C \subset \mathfrak{t}^{*}$ the corresponding Weyl chamber.  The dual object
$\widehat{G}$ of $G$ is in one-to-one correspondence with $\overline{C} \cap I^{*}$. (When $G$ is semisimple and simply connected,
there exist integral forms $\rho_1,...,\rho_l \in \mathfrak{t}^{*}$ satisfying
$2\langle \rho_j,\alpha_k \rangle / \langle \alpha_k,\alpha_k \rangle = \delta_{jk}$ for $j,k \in \{1,...,l\}$ and
$\overline{C} \cap I^{*} = \{n_1\rho_1 + ... + n_l\rho_l: \text{ each } n_j \in \mathbb{Z} \text{ and } n_j \geq 0\}$.)

Fix $\omega \in C \cap I^{*}$ and for $N = 1,2,3,...$ set $\Lambda_{N}=\{\lambda\in\overline{C}\cap I^{*}: \lambda \leq N\omega\}$.
The associated polyhedral Dirichlet kernel on $G$ is $\boldsymbol{D}_N = \sum_{\lambda\in\Lambda_N}d_{\lambda}\chi_{\lambda}$
where $N=1,2,3,...$, and the convolution products $S_{N}f = \boldsymbol{D}_N \ast f$ define the sequence of associated
polyhedral partial sums for the Fourier series of $f \in L^2(G)$.

Using the Weyl integration formula on $G$ \cite{BD}, the $N$th Fourier polyhedral partial sum of $f$ at
$x \in G$ is given by
\begin{align*}
S_{N} & f(x) =(\boldsymbol{D}_{N}\ast f)(x)\\
 & =\intop_{G}\boldsymbol{D}_{N}(y)f(y^{-1}x)d\mu(y)\\
 & =\frac{1}{|W|}\intop_{T}\det(E_{G/T}-\textrm{Ad}_{G/T}(\theta^{-1}))
    \intop_{G}\boldsymbol{D}_{N}(y{\theta}y^{-1})f(y\theta^{-1}y^{-1}x)d\mu(y)d\theta\\
 & =\frac{1}{|W|}\intop_{T}\eta(\theta)\boldsymbol{D}_{N}(\theta)[Q_{x}f](\theta)d\theta
\end{align*}
where $\eta(\theta) = \det(E_{G/T}-\textrm{Ad}_{G/T}(\theta^{-1}))$ and $[Q_{x}f](\theta)=\intop_{G}f(y\theta^{-1}y^{-1}x)d\mu(y)$.

Let $\pi(\lambda)$ denote the continuous, irreducible, unitary representation of $G$ corresponding to
$\lambda \in \overline{C}\cap I^{*}$, let $H_{\pi(\lambda)}$ be the finite dimensional subspace of $L^{2}(G)$ generated by
the coordinate functions of $\pi(\lambda)$ \cite[p.24]{HR}, and let $P_{\lambda}f = d_{\lambda}(\chi_{\lambda}\ast f)$
be the orthogonal projection of $L^{2}(G)$ onto $H_{\pi(\lambda)}$.  Then
\begin{equation*}
S_{N}f(x) =\sum_{\lambda \in \Lambda_{N}}(P_{\lambda}f)(x) = \sum_{j=0}^{N}(\Gamma_{j}f)(x)
\end{equation*}
with $\Gamma_{n}f = \sum P_{\lambda}f$ where the sum is over all $\lambda \in \overline{C}\cap I^{*}$ such that
$(n-1)\omega < \lambda \leq n\omega$.  Let $\delta_{h}f(x)= f(x)-f(h^{-1}x)$ be the difference operator
on $L^{2}(G)$ and let the integral modulus of continuity of $f \in L^{2}(G)$ be given by
\begin{equation*}
  \Omega(f,t)=\sup\{\Vert\delta_{h}f\Vert_{L^{2}(G)}: h = \exp(X), X \in \mathfrak{g}, \textrm{ and }0 < \Vert X \Vert \leq t\};
\end{equation*}
here $\Vert \cdot \Vert$ denotes the norm on $\mathfrak{g}$ generated by the inner product $\langle \cdot,\cdot \rangle$.

An alternate theory for almost everywhere convergence of Fourier series of $f \in L^2(G)$ can be based on spherical partial
sums (cf. \cite{CK} and \cite{GT2}).  Let
\begin{equation*}
  \rho = \frac{1}{2}\sum_{\alpha \in P_{+}}\alpha
\end{equation*}
and for $N=1,2,3,...$ define $\widetilde{\Lambda}_{N}=\{\lambda\in\overline{C}\cap I^{*}: \Vert \lambda-\rho \Vert \leq N\}$,
the spherical Dirichlet kernel $\widetilde{\boldsymbol{D}}_N = \sum_{\lambda\in\widetilde{\Lambda}_N}d_{\lambda}\chi_{\lambda}$,
and the spherical partial sums
\begin{equation*}
\widetilde{S}_{N}f = \widetilde{\boldsymbol{D}}_{N}\ast f = \sum_{\lambda \in \widetilde{\Lambda}_{N}}(P_{\lambda}f)
    = \sum_{j=0}^{N}(\widetilde{\Gamma}_{j}f)
\end{equation*}
with $\widetilde{\Gamma}_{n}f = \sum P_{\lambda}f$ where the sum is over all $\lambda \in \overline{C}\cap I^{*}$ such that
$n-1 < \Vert \lambda-\rho \Vert \leq n$.  We see no essential difference between the behavior of polyhedral and spherical partial
sums in the results of this work.

\section{Almost everywhere convergence results}

\begin{theorem}
Let $G$ be a compact, connected Lie group and $f\in L^{2}(G)$. If
\begin{equation*}
  \intop_{0}^{1}\frac{\Omega^{2}(f,t)}{t}dt < \infty,
\end{equation*}
then $S_{N}f(x)\rightarrow f(x)$ and $\widetilde{S}_{N}f(x)\rightarrow f(x)$ for almost every $x \in G$.
\end{theorem}

\begin{proof}
We give the proof for polyhedral Fourier partial sums of $f$; the case of spherical Fourier partial sums is completely
analogous.  Note that by Jackson's theorem for compact, connected Lie groups \cite{CK}, there is a constant $B > 0$ such that
\begin{align*}
\intop_{0}^{1}\frac{\Omega^{2}(f,t)}{t}dt & =\sum_{k=0}^{\infty}\intop_{2^{-(k+1)}}^{2^{-k}}\frac{\Omega^{2}(f,t)}{t}dt\\
 & \geq \log(2)\sum_{k=1}^{\infty}\Omega^{2}(f,2^{-k})\\
 & \geq B\log(2)\sum_{k=1}^{\infty}(E_{2^{k}}^{2}f)
\end{align*}
where
\begin{equation*}
  E_{M}(f)=\inf\{\Vert f-P \Vert_{L^{2}(G)}: P \in\bigoplus_{\lambda\in\Lambda_M}H_{\pi(\lambda)}\textrm{ }\}
\end{equation*}
is the best approximation of $f$ in $\bigoplus_{\lambda\in\Lambda_M}H_{\pi(\lambda)}$.  Since
\begin{align*}
E_{2^{k}}^{2}(f) & =\Vert f-S_{2^{k}}f\Vert_{L^{2}(G)}^{2}\\
 & =\sum_{j=2^{k}+1}^{\infty}\Vert\varGamma_{j}f\Vert_{L^{2}(G)}^{2},
\end{align*}
it follows that
\begin{align*}
\sum_{k=1}^{\infty}E_{2^{k}}^{2}(f) & =\sum_{j=2}^{\infty}\left\lfloor \log_{2}(j)\right\rfloor \Vert\varGamma_{j}f\Vert_{L^{2}(G)}^{2}\\
 & \geq\frac{\log_{2}(e)}{2}\sum_{j=2}^{\infty}\log(j)\Vert\varGamma_{j}f\Vert_{L^{2}(G)}^{2},
\end{align*}
and therefore
\begin{equation*}
  \infty > \intop_{0}^{1}\frac{\Omega^{2}(f,t)}{t}dt \geq \frac{B}{2}\sum_{j=2}^{\infty}\log(j)\Vert\varGamma_{j}f\Vert_{L^{2}(G)}^{2}.
\end{equation*}
Consequently, a general version of the Rademacher-Menshov theorem \cite{M} implies
\begin{equation*}
  \lim_{N \rightarrow \infty}S_{N}f(x)=\lim_{N \rightarrow \infty}\sum_{j=0}^{N}(\varGamma_{j}f)(x)=f(x)
\end{equation*}
for almost every $x \in G$.
\end{proof}

It should be noted that in \cite{D}, Dai established almost everywhere convergence results for Fourier-Laplace series
of functions on spheres which are analogous to the preceding theorem.

\begin{cor}
Let $f\in L^{2}(G)$. If $\Omega(f,t)= O(t^{\alpha})$ for some $\alpha > 0$, then $S_{N}f(x)\rightarrow f(x)$
and $\widetilde{S}_{N}f(x)\rightarrow f(x)$ for almost every $x\in G.$
\end{cor}

\begin{cor}
If $f\in\textrm{Lip}_{\alpha}(G)$ for some $\alpha > 0$, then $S_{N}f(x)\rightarrow f(x)$ and
$\widetilde{S}_{N}f(x)\rightarrow f(x)$ for almost every $x \in G$.
\end{cor}

\section{Divergence of Fourier partial sums on a countable subset}

To begin this section, we review notation and results of \cite[p.84f]{BD} and \cite[p.125f]{HR}.
Equip the two-dimensional special unitary group $SU(2)$ with the left and right translation invariant metric $d$ given by
\begin{equation*}
  d(x,y)= \sqrt{\frac{1}{2}\textrm{tr}\big((x-y)(x-y)^{*}\big)}.
\end{equation*}
Let $0 < \alpha < 1$ and let $f$ be a real function on $SU(2)$. If there exists a number $M > 0$ such that
\begin{equation*}
\vert f(x)-f(y) \vert \leq Md^{\alpha}(x,y)
\end{equation*}
for all $x,y \in SU(2)$, then $f$ is an $\alpha$-H\"{o}lder continuous function on $SU(2)$ and we write
$f\in\textrm{Lip}_{\alpha}(SU(2))$. A real function $f$ on $SU(2)$ is central if, for $\mu$-almost
every $x \in SU(2)$, $f(yxy^{-1})=f(x)$ for all $y \in SU(2)$.  In particular, since every $x \in SU(2)$ is
diagonalizable via a similarity transformation:
\begin{equation*}
yxy^{-1} =
  \begin{bmatrix}
    &e^{i\theta} &0 \\
    &0  &e^{-i\theta}
  \end{bmatrix}
  \equiv \omega(\theta),
\end{equation*}
where $y \in SU(2)$ and $e^{\pm i\theta}$ are the eigenvalues of $x$, it follows that if $f$ is central then for
$\mu$-almost every $x \in SU(2)$,
\begin{equation*}
f(x) = f(\omega(\theta))
\end{equation*}
with $\theta \in [0,\pi]$. With this notation, the Weyl integral formula is explicitly
\begin{equation}\label{CentralIntegral}
  \int_{SU(2)}g(x)d\mu(x)=\frac{2}{\pi}\int_{0}^{\pi}g(\omega(\theta))\sin^2(\theta)d\theta
\end{equation}
when $g$ is a central function in $L^2(SU(2))$.

Let $\{\pi_n\}_{n=0}^{\infty}$ denote the family of all (inequivalent) continuous, irreducible, unitary representations
of $SU(2)$; here $\pi_n$ has dimension $n+1$ and its character $\chi_n = \textrm{trace}(\pi_n)$ is the continuous central
function on $SU(2)$ given by
\begin{equation*}
\chi_n(x) = \chi_n(\omega(\theta)) = \frac{\sin((n+1)\theta)}{\sin(\theta)}
\end{equation*}
where $e^{\pm i\theta}$ are the eigenvalues of $x$.  It follows that the Dirichlet kernel
$\{\mathbf{D}_N\}_{N=0}^{\infty}$ on $SU(2)$ is the sequence of continuous central functions given by
\begin{equation}\label{Dirichlet}
\mathbf{D}_N(x)= \mathbf{D}_N(\omega(\theta)) = \sum_{n=0}^N (n+1)\chi_n(\omega(\theta)) =
    \frac{-1}{2\sin(\theta)}D_{N+1}^{\prime}(\theta)
\end{equation}
where
\begin{equation*}
  D_n(t)=1+2\sum_{j=1}^{n}\cos(jt)=\frac{\sin((2n+1)t/2)}{\sin(t/2)}
\end{equation*}
is the Dirichlet kernel on $[-\pi,\pi]$.

\begin{theorem}\label{divergence}
Let $\alpha\in(0,1)$ and let $\{x_{i}\}_{i=1}^{\infty}$ be any countable subset of $SU(2)$. Then there exists a function
$f\in\textrm{Lip}_{\alpha}(SU(2))$ such that
\[
\underset{N\geq1}{\sup}|S_{N}f(x_{i})|=\infty
\]
for all $i=1,2,3,\ldots$.
\end{theorem}

\begin{proof}
Observe that for $\alpha \in (0,1)$, $\textrm{Lip}_{\alpha}(SU(2))$ is a Banach space with norm
\begin{equation*}
\Vert f\Vert_{\textrm{Lip}_{\alpha}(SU(2))}=\underset{x\in SU(2)}{\sup}|f(x)|+\underset{\underset{x\neq y}{x,y\in}SU(2)}{\sup}\frac{|f(x)-f(y)|}{d^{\alpha}(x,y)}.
\end{equation*}
Fix $x \in SU(2)$ and $n \in \mathbb{N}$, and set $\Phi_{n}^{x}(f)=S_{n}f(x)$.  Each $\Phi_{n}^{x}$ is a
bounded linear functional on $\textrm{Lip}_{\alpha}(SU(2))$ of norm
\begin{align*}
  \Vert\Phi_{n}^{x}\Vert &
    =\sup\left\{ |S_{n}f(x)|:f\in\textrm{Lip}_{\alpha}(SU(2)),\;\Vert f\Vert_{\textrm{Lip}_{\alpha}(SU(2))}\leq1\right\} \\
  &\leq \Vert\boldsymbol{D}_{n}\Vert_{L^{1}(SU(2))}.
\end{align*}
Specializing to the case when $x = e$, the identity matrix in $SU(2)$, and $f \in L^2(SU(2))$ is central, we have
\begin{align*}
  \Phi_{n}^{e}(f)=(f\ast\boldsymbol{D}_{n})(e)& =
      \int_{SU(2)}f(y)\boldsymbol{D}_{n}(y)d\mu(y) \\
      & =\frac{2}{\pi}\int_{0}^{\pi}f(\omega(\theta))\boldsymbol{D}_{n}(\omega(\theta))\sin^2(\theta)d\theta
\end{align*}
by (\ref{CentralIntegral}).
It then follows from (\ref{Dirichlet}) that
\begin{align*}
\Phi_{n}^{e}(f)
  =&\frac{1}{\pi}\int_{0}^{\pi}f(\omega(\theta))\cos^{2}\left(\frac{\theta}{2}\right)D_{n+1}(\theta)d\theta\\
   &\quad - \frac{(2n+3)}{\pi}\int_{0}^{\pi}f(\omega(\theta))\cos\left(\left(n+\frac{3}{2}\right)\theta\right)\cos\left(\frac{\theta}{2}\right)d\theta.
\end{align*}

For $n \geq 2$, the absolute maxima and minima of the function $h_n(\theta)=\cos\left(\left(n+\frac{3}{2}\right)\theta\right)$
on $[0,\pi]$ occur at the endpoints of the intervals
\begin{equation*}
  I_{k}=\left[\frac{2k\pi}{2n+3},\frac{2(k+1)\pi}{2n+3}\right]
\end{equation*}
where $k \in \{0,1,2,...,n\}$.  Let $g_n$ be the sawtooth function on $[0,\pi]$ determined by
$g_n\left(\frac{2k\pi}{2n+3}\right) = (-1)^k$ for $0\leq k \leq n+1$, $g_n(\pi)=0$, and $g_n$ is piecewise linear
between these points.  Define a central function $f_n$ on $SU(2)$ by $f_n(\omega(\theta)) = g_n(\theta)$ for
$\theta \in [0,\pi]$. It is easy to see that each $f_n$ belongs to $\textrm{Lip}_{\alpha}(SU(2))$; in fact,
\begin{equation*}
  \frac{\vert f_n(x)-f_n(y) \vert}{d^{\alpha}(x,y)} \leq \left(\frac{\pi}{2}\right)^{\alpha}\left(\frac{2\pi}{2n+3}\right)^{1-\alpha}
    \leq \pi
\end{equation*}
for all distinct matrices $x$ and $y$ in $SU(2)$.

Since $g_n(\theta)\cos\left(\left(n+\frac{3}{2}\right)\theta\right) \geq g_n^2(\theta) \geq 0$ on each interval $I_k$ and on
$\left[\frac{2(n+1)\pi}{2n+3},\pi\right]$, and since the function $\theta \mapsto \cos(\theta/2)$ is positive and decreasing
on $[0,\pi)$,
\begin{align*}
   \int_{I_k}f_n(\omega(\theta))\cos\left(\left(n+\frac{3}{2}\right)\theta\right)\cos(\theta/2)d\theta
        &\geq \int_{I_k}g_n^2(\theta)\cos(\theta/2)d\theta \\
        &\geq \cos\left(\frac{(k+1)\pi}{2n+3}\right)\int_{I_k}g_n^2(\theta)d\theta \\
        &= \frac{2\pi}{3(2n+3)}\cos\left(\frac{(k+1)\pi}{2n+3}\right)
\end{align*}
for all $k \in \{0,1,2,...,n\}$ and
\begin{equation*}
\int_{\frac{2(n+1)\pi}{2n+3}}^{\pi}f_n(\omega(\theta))\cos\left(\left(n+\frac{3}{2}\right)\theta\right)\cos(\theta/2)d\theta \geq 0.
\end{equation*}
Adding these inequalities we obtain
\begin{align*}
\intop_{0}^{\pi}f_n(\omega(\theta))\cos\left(\left(n+\frac{3}{2}\right)\theta\right)\cos\left(\frac{\theta}{2}\right)d\theta
    & \geq \frac{2}{3}\left(\frac{\pi}{2n+3}\right)\sum_{k=1}^{n+1}\cos\left(\frac{k\pi}{2n+3}\right)\\
    & = \frac{2}{3}\left(\frac{\pi}{2n+3}\right)\left\{ D_{n+1}\left(\frac{\pi}{2n+3}\right)-1\right\},
\end{align*}
and hence
\begin{equation*}
\biggl|\frac{(2n+3)}{\pi}\int_{0}^{\pi}f_n(\omega(\theta))\cos\left(\left(n+\frac{3}{2}\right)\theta\right)\cos\left(\frac{\theta}{2}\right)d\theta\biggr|
\geq\frac{2}{3}\left\{ D_{n+1}\left(\frac{\pi}{2n+3}\right)-1\right\} .
\end{equation*}
Because the function $\theta\mapsto f_n(\omega(\theta))\cos^{2}(\theta/2)$ is uniformly bounded by $1$ on $[0,\pi]$,
\begin{align*}
\biggl |\frac{1}{\pi}\int_{0}^{\pi}f_n(\omega(\theta))\cos^{2}(\theta/2)D_{n+1}(\theta)d\theta \biggr|
    & \leq \frac{1}{\pi}\intop_{0}^{\pi}|D_{n+1}(\theta)|d\theta \\
    & = \frac{4}{\pi^{2}}\log(n+1)+o(1)
\end{align*}
as $n\rightarrow\infty$. Consequently
\begin{equation*}
\frac{|\Phi_{n}^{e}(f_n)|}{\Vert f_n\Vert_{\textrm{Lip}_{\alpha}(SU(2))}}
    \geq \frac{\frac{2}{3}\left\{ D_{n+1}\left(\frac{\pi}{2n+3}\right)-1\right\}
        -\left(\frac{4}{\pi^{2}}\log(n+1)+o(1)\right)}{1+\pi}.
\end{equation*}
But $D_{n+1}\left(\frac{\pi}{2n+3}\right) = \left(\sin\left(\frac{\pi}{2(2n+3)}\right)\right)^{-1} \geq \frac{2(2n+3)}{\pi}$
and hence
\begin{equation*}
\Vert\Phi_{n}^{e}\Vert =
\sup\left\{\frac{|\Phi_{n}^{e}(f)|}{\Vert f\Vert_{\textrm{Lip}_{\alpha}(SU(2))}}:\quad f\in\textrm{Lip}_{\alpha}(SU(2)),
    \quad f\neq 0 \right\}
\end{equation*}
is asymptotically bounded below by
\begin{equation*}
    \frac{\frac{2}{\pi}(2n+3)-\frac{4}{\pi^{2}}\log(n+1)}{1+\pi}
\end{equation*}
as $n\rightarrow\infty$.  Thus the sequence of bounded linear functionals
\begin{equation*}
    \Phi_{n}^{e}(f)=S_{n}f(e)
\end{equation*}
is not uniformly bounded on the Banach space $\textrm{Lip}_{\alpha}(SU(2))$ as $n\rightarrow\infty$. By the uniform
boundedness principle
\begin{equation*}
    \underset{n\geq1}{\sup}|S_{n}f(e)| = \infty
\end{equation*}
for all $f$ belonging to some dense $G_{\delta}$ set in $\textrm{Lip}_{\alpha}(SU(2))$.

If $z \in SU(2)$, define the left translation operator $L_z$ on $\textrm{Lip}_{\alpha}(SU(2))$ by $L_{z}f(y) = f(zy)$
for all $y \in SU(2)$.  For each element of the countable subset $\{x_{i}\}_{i=1}^{\infty}$ of $SU(2)$ and each
$n \geq 1$, observe that
\begin{equation*}
\frac{|\Phi_{n}^{x_i}\left(L_{{x_i}^{-1}}f_{n}\right)|}{\Vert L_{{x_i}^{-1}}f_{n}\Vert_{\textrm{Lip}_{\alpha}(SU(2))}}
  =\frac{|\Phi_{n}^{e}\left(f_{n}\right)|}{\Vert f_{n}\Vert_{\textrm{Lip}_{\alpha}(SU(2))}},
\end{equation*}
so there corresponds a dense $G_{\delta}$ subset $E_{x_i}$ of $\textrm{Lip}_{\alpha}(SU(2))$ such that
\begin{equation*}
    \underset{n\geq1}{\sup}|S_{n}f(x_i)|=\infty
\end{equation*}
for all $f \in E_{x_i}$. By the Baire category theorem $E=\bigcap_{i=1}^{\infty}E_{x_{i}}$ is dense in
$\textrm{Lip}_{\alpha}(SU(2))$. In particular, $E$ is nonempty and any $f \in E$ gives the desired conclusion.
\end{proof}

A general theorem of Belen'kii \cite{B} guarantees that if $F:[-1,1] \rightarrow \mathbb{R}$ satisfies a
Dini-Lipschitz condition and the Fourier-Jacobi partial sums
\begin{equation*}
  \{s_N^{(\alpha,\beta)}(F;\pm 1)\}_{N=0}^{\infty}
\end{equation*}
at $\pm 1$ converge for some $\alpha > -1$ and $\beta > -1$, then the corresponding Fourier-Jacobi series of $F$
converges uniformly to $F$ on $[-1,1]$.  This leads directly to the following result on $SU(2)$.
\begin{theorem}\label{central}
  Let $f$ be a central function in $\textrm{Lip}_{\alpha}(SU(2))$ for some $\alpha \in (0,1)$.  Then:
  \begin{enumerate}[label=(\alph*)]
    \item $S_Nf(x) \rightarrow f(x)$ uniformly outside any open set containing $\{e,-e\}$;
    \item $S_Nf(x) \rightarrow f(x)$ uniformly on $SU(2)$ if $\{S_Nf(\pm e)\}_{N=0}^{\infty}$ converge.
  \end{enumerate}
\end{theorem}
This suggests the question: Can the word ``central'' be deleted from the hypothesis of Theorem
\ref{central} and still obtain conclusions (a) and (b)?  Theorem \ref{divergence} shows that there is no
possibility of such an analogue of Theorem \ref{central} for general functions in $\textrm{Lip}_{\alpha}(SU(2))$ for
some $\alpha \in (0,1)$.  The points of divergence for the Fourier partial sums of such a noncentral
function need no longer be at the ``poles'' $\pm e$ of the ``sphere'' $SU(2)$. According to Theorem \ref{divergence},
points of divergence for $\textrm{Lip}_{\alpha}(SU(2))$ functions can be dense in $SU(2)$.

\end{document}